\documentclass[leqno,12pt]{article}
\newcommand{\bas}{\begin{eqnarray*}} \newcommand{\eas}{\end{eqnarray*}}
\def\tto{\;{\lower 1pt \hbox{$\rightarrow$}}\kern -12pt
           \hbox{\raise 2.8pt \hbox{$\rightarrow$}}\;}
\def\for{\hskip0.9pt|\hskip0.9pt}
\newcommand{\eps}{\epsilon}

\newcommand{\bd}{\begin{displaymath}} \newcommand{\ed}{\end{displaymath}}
\newcommand{\be}{\begin{equation}} \newcommand{\ee}{\end{equation}}

\def\bx{\bar x} \def\by{\bar y}  \def\bp{\bar p}
   \def\bp{\bar p}
 
\def\gph{\mathop{\rm gph}\nolimits}

\def\reg{\mathop{\rm mod}\nolimits} 
\def\calm{\mathop{\rm clm}\nolimits}
 \def\Int{\mathop{\rm int}\nolimits}

\def\plus{{\scriptscriptstyle +}} \def\minus{{\scriptscriptstyle -}}
\def\starplus{{*\hskip-0.6pt{\raise0.5pt\hbox{$\plus$}}}}
\def\starminus{{*\hskip-0.6pt{\raise0.5pt\hbox{$\minus$}}}}

\def\epsilon{\varepsilon}             \def\phi{\varphi}
\def\reals{{I\kern-.35em R}}

\def\iff{~~\Leftrightarrow~~} 
\def\text#1{\;\,\hbox{#1}\;\,}       
\def\Text#1{\quad\;\,\hbox{#1}\quad\;\,}
\def\state #1. { \noindent{\bf#1.\enspace}}
\def\lset{\big\{\,}       \def\rset{\big\}}
       
\outer\def\proclaim #1. #2\par{\medbreak \noindent{\bf#1.\enspace}{\sl#2}\par
  \ifdim\lastskip<\medskipamount \removelastskip\penalty55\medskip\fi}
\def\qed{\hfill{$\vcenter{\hrule height1pt \hbox{\vrule width1pt height5pt
   \kern5pt \vrule width1pt} \hrule height1pt}$} \medskip}
\def\low#1{{\lower1pt \hbox{$\scriptstyle #1$}}}
\def\downto{{\raise 1pt \hbox{$\scriptstyle \,\searrow\,$}}}

\def\bx{\bar x} \def\by{\bar y}

\def\gph{\mathop{\rm gph}\nolimits} \def\lip{\mathop{\rm lip}\nolimits}

\def\reg{\mathop{\rm reg}\nolimits}

\def\plus{{\scriptscriptstyle +}} \def\minus{{\scriptscriptstyle -}}

\def\epsilon{\varepsilon}             \def\phi{\varphi}
\def\reals{{I\kern-.35em R}}

\def\iff{~~\Leftrightarrow~~} 
\def\text#1{\;\,\hbox{#1}\;\,}       
\def\Text#1{\quad\;\,\hbox{#1}\quad\;\,}
\def\state #1. { \noindent{\bf#1.\enspace}}
\def\lset{\big\{\,}       \def\rset{\big\}}
       
\outer\def\proclaim #1. #2\par{\medbreak \noindent{\bf#1.\enspace}{\sl#2}\par
  \ifdim\lastskip<\medskipamount
  \removelastskip\penalty55\medskip\fi}
\def\qed{\hfill{$\vcenter{\hrule height1pt \hbox{\vrule width1pt height5pt
   \kern5pt \vrule width1pt} \hrule height1pt}$} \medskip}
\def\low#1{{\lower1pt \hbox{$\scriptstyle #1$}}}
\def\downto{{\raise 1pt \hbox{$\scriptstyle \,\searrow\,$}}}
\def\tto{\;{\lower 1pt \hbox{$\rightarrow$}}\kern -12pt
           \hbox{\raise 2.8pt \hbox{$\rightarrow$}}\;}
\def\for{\hskip0.9pt|\hskip0.9pt}

\def\ball{{I\kern -.35em B}}

\begin{document}

\title{ An Inverse Function Theorem \\for Metrically Regular Mappings}

\author{A. L. Dontchev  \\
Mathematical Reviews  \\
 Ann Arbor, MI 48107-8604,  \\
{\tt ald@ams.org} }

\date{August 22, 2002}

\maketitle

\begin{abstract}
We prove that if a mapping $F:X \tto Y$, where $X$ and $Y$
are  Banach spaces,  is metrically
regular at $\bx$ for $\by$ and its inverse $F^{-1}$ is convex and closed valued
locally around $(\bx, \by)$, then  for any
function $G:X\to Y$ with $\lip G(\bx)\cdot \reg F(\bx\for\by)) < 1$, the mapping
$(F+G)^{-1}$ has a continuous local selection $x(\cdot)$
around $(\bx, \by+G(\bx))$ which is also calm. 
 \end{abstract}
\bigskip\medskip

{\small \noindent {\bf Key words:} set-valued mapping, 
metric regularity, 
continuous selections, inverse/implicit function theorem.    

}

\newpage

\section{ Introduction }

The classical inverse function theorem  stated for a function
 $f:X  \to Y$, with $X$ and $Y$ Banach spaces, 
 assumes
that  $f$ is continuously
differentiable in a neighborhood of a given reference point $\bx$
and, most importantly, the Fr\'echet
derivative $ \nabla f(\bx)$ has a linear
and bounded inverse; then the theorem claims that
 there exist neighborhoods $U$ of $\bx$ and $V$ of $\by:=f(\bx)$
such that the mapping
$$V \ni y \mapsto f^{-1}(y)\cap U \leqno{(1.1)}$$
is single valued (a function defined on $V$) which is moreover
continuously differentiable (${\cal C}^1$)  in $V$ and whose derivative
is the inverse of $\nabla f$. It is perhaps less known
that an inverse function type  theorem may be obtained
when  the Jacobian $ \nabla f(\bx)$ is merely surjective.
Indeed, in this case the mapping (1.1), although  in general
set-valued, may have a local single-valued selection $x(\cdot)$, that is, 
a function $x(\cdot)$ exists with $x(y) \in f^{-1}(y)
\cap U$ for all $y \in V,$ which  is continuously differentiable in $V$.
In other words, a smooth  inverse function exists but it is only a part
of the inverse $f^{-1}$ which may be set-valued.
 The precise result is as follows:

\proclaim  Theorem 1.1. Let $X$ and $Y$ be Hilbert
spaces and let $f:X \to Y$ be a function which is ${\cal C}^1$
around $\bx$ and such that the derivative
$B:=\nabla f(\bx)$ is surjective. Then there exist a neighborhood
$V$ of $\by:= f(\bx)$ and a ${\cal C}^1$
 function $x:V \to X$  such that 
$$x(\by) = \bx \quad \Text{ and } \quad f(x(y)) = y \text{ for every}   y \in V,$$
 and moreover $ \nabla x(\by) = (B^*B)^{-1}B^*$.

{\bf Proof.} In terms of the adjoint operator $B^*$ consider the mapping
$$(x,u) \mapsto g(x,u): = \left( \begin{array}{ll}  
x+ B^*u \\ f(x) \end{array}\right),$$
which satisfies $g(\bx, 0) = (\bx,\by)$ and whose Jacobian is
$$J = \left( \begin{array}{ll} I & B^*\\B& 0 \end{array}\right). $$
It is well known that, in Hilbert spaces, when $B$ is surjective than
the operator $J$ is invertible in the sense that $J^{-1}$ is  linear and bounded
 from $X\times Y$ into itself. Hence, by the classical inverse function theorem
 the mapping $g^{-1}$, when restricted to
a neighborhood of the point $((\bx, 0), (\bx, \by))$ in its graph,
is single-valued and continuously differentiable. In particular, for some
neighborhoods $U$ of $\bx$ and $V$ of $\by$, the function $x(y) = \xi(\bx,y)$ satisfies
$y = f(x(y))\cap U$
for $y \in V$. It remains to observe that $B^*B$ is invertible and, from the 
equation $B^*f(x(y)) = B^* y$, 
the derivative of $x(\cdot)$ with respect to $y$ 
satisfies $B^*B\nabla x(\by) = B^*$. \qed

The implicit function theorem corresponding to Theorem 1.1 
is easy to  prove
by using the standard passage from inverse to implicit function theorems. 
An interesting reading about the history and theory of 
various implicit function theorems
is the recent book \cite{KP}.

If $X$ and $Y$ are arbitrary Banach spaces, it is in
our opinion  quite unlikely
that a result of the form of  Theorem~1.1 holds; however, we do not know
a counterexample. Still,  in Banach spaces the surjectivity of
the Jacobian implies the existence of a selection of (1.1) which 
may be not smooth but it is continuous and calm.
Specifically, we have 

\proclaim Theorem  1.2. Let  $X$ and $Y$ be Banach spaces
 and let $f:X \to Y$ be a function which is strictly differentiable
at $\bx$ and such that the strict derivative
$\nabla f(\bx)$ is surjective. Then there exist a neighborhood
$V$ of $\by:= f(\bx)$,  a continuous function $x:V \to X$
and a constant $\gamma > 0$   such that 
$$ f(x(y)) = y  \Text{and} \|x(y) - \bx\| \leq \gamma\|y - \by\| \,
\text{  for every}   y \in V. $$

Theorem 2.1 was communicated to the author by  H. Sussmann \cite{S} who put it in the
context of the Lyusternik theorem. A  version of this theorem, without 
the estimate,  appeared in  \cite{G}.

In this paper we will obtain Theorem 1.2 as a corollary of the following 
more general result:
Let a  set-valued mapping
$F:X \tto Y$ be metrically
regular at $\bx$ for $\by$ and its inverse $F^{-1}$ be convex and closed valued
locally around $(\bx, \by)$. Then $F^{-1}$ has a continuous local selection $x(\cdot)$
around $(\bx, \by)$ which  is calm.
Moreover,  for any 
function $G:X\to Y$ with $\lip G(\bx)\cdot \reg F(\bx\for\by)) < 1$, the mapping
$(F+G)^{-1}$ has a continuous local selection $x(\cdot)$
around $(\bx, \by+G(\bx))$ which is calm. Here $\reg F(\bx\for\by)$
is the modulus of strong regularity of $F$ which is defined in further lines
and $\lip G(\bx)$ is the Lipschitz modulus of $G$, see Section 3 for a definition.

The general paradigm behind our result is the same as in standard 
inverse/implicit mapping theorems: Suppose a mapping $f$ can be
represented as the sum $f=F+G$, where $F$ is ``nice"
(metrically regular with  locally convex and closed inverse) so that
it has an ``inverse" and $G$ is ``small"
(small Lipschitz constant); then $f$ has an ``inverse" as well. 
In Theorem 1.2 $F$ is the Jacobian mapping
of the function $f$ and $G$ is the difference $f-F$,
as in the classical inverse function theorem.

In the remaining part of this section we describe the notation 
and terminology we use, which is consistent with
the book \cite{RW}, and briefly discuss some related  results.
Throughout, unless stated otherwise,
 $X$ and $Y$ are real Banach spaces with norms $\|\cdot \|$
and closed unit balls $\ball$; a ball centered at $a$ with radius $r$
is  $\ball_r(a)$. The distance from a point $x$ to a
set $A$ is denoted by $d(x, A)$. The notation
$F: X \tto Y$ means that $F$ is a {\em set-valued} 
mapping from $X$ to the subsets of $Y$;
if $F$ is a function, that is, for each $x\in X$
the set of values $F(x)$ consists of no more than one
element, then we write $F: X \to Y$. The {\em graph } of 
$F$ is  $\gph F = \{ (x,y) \mid y \in F(x) \}$ and its {\em inverse}
$F^{-1}$ is defined as $x \in F^{-1}(y) \iff y \in F(x)$. 
A mapping $F:X \tto Y$ with $(\bx, \by) \in \gph F$
has a {\em local selection around} $(\bx, \by)$
 if there exist neighborhoods $U$ of $\bx$
and $V$ of $\by$ and a function $s:U \to V$
such that $s(\bx) = \by$ and $s(x) \in F(x)\cap V$ for all $x \in U$.
A mapping $F$ has a {\em continuous local selection}
 if it has a selection $s:U\to V$  which
is continuous in $U$.
  We say that a function
$f$ is  ${\cal C}^1$  when $f$ is continuously
differentiable (smooth). The (strict or  Fr\'echet) derivative of $f$ at $\bx$
is denoted by $\nabla f(\bx)$. A function $f:X \to Y$ is 
{\em calm at} $\bx$ when there exist a neighborhood
$V$ of $\bx$ and a constant $\gamma > 0$   such that 
$$ \|f(x) - f(\bx)\| \leq \gamma\|x - \bx\| \,
\text{  for every}   x \in V. \leqno{(1.2)}$$
 The infimum of $\gamma$ for which (1.2) holds is
called {\em modulus of calmness} and is denoted by $\calm f(\bx)$.

A mapping $F: X \tto Y$ is said to be {\em metrically
regular } at $\bx$ for $\by$ if there exists a constant $\kappa > 0$ such that
$$d(x, F^{-1}(y)) \leq \kappa d(y, F(x))
\text{for all $(x,y)$ close to $(\bx,\by).$}     \leqno(1.3)
$$
The infimum of $\kappa$ for which (1.3) holds is the {\em modulus
of metric regularity} which we denote by $\reg F(\bx\for\by)$.

The concept of metric regularity  
has its roots in the work of L. A. Lyusternik in the 30s 
and L. Graves in the 50s. It has been playing a central role in 
optimization for obtaining necessary  conditions for extremum.
 In a very general setting, the metric regularity  is
``stable under linearization,'' in the line of the Lyusternik-Graves
theorem and its various extensions. Further, it is ``robust"; that is,
if it holds for a mapping at certain point, it holds 
also in  ``neighborhoods" of the mapping and the  point.   As shown recently in
\cite{rad}, in a finite-dimensional setting,
the distance from a given metrically regular mapping to
 the set of mappings that are not
metrically regular, measured by the Lipschitz modulus of the perturbation, is
equal to the reciprocal to the modulus of metric regularity.
A discussion of  various developments around
this concept  has recently 
been given  by Ioffe \cite{Ioffe}, details 
 are also available in \cite{RW}.  

In finite dimensions,  more can be said about 
functions implicitly defined by metrically regular mappings.
If $f:\reals^n \to \reals^n$ is continuously differentiable around $\bx$,
then the metric regularity of $f$ at $\bx$  simply means
that the Jacobian $\nabla f(\bx)$ is a nonsingular matrix and
then the localization of $f^{-1}$ around the point $(f(\bx), \bx)$
is single-valued and ${\cal C}^1$. 
The equivalence of the metric regularity with the {\em Lipschitz continuous}
single-valued localization of $f^{-1}$ is actually valid for  more
general set-valued mappings of the form $f + N_C$ where $f$ is a smooth function
and $N_C$ is the normal cone mapping to a convex polyhedral set $C$.
This inverse function theorem for variational inequalities
was established in \cite{AT} together with a formula for the 
Lipschitz modulus of the localization. Here the theory of inverse function
for metrically regular mappings  merges with another 
fundamental result, due to S. Robinson \cite{Rob},
regarding the ``stability under linearization'' of the property of existence of 
a Lipschitz continuous single-valued localization.
The existence of a single-valued
and continuous localization of the inverse also holds  for merely
continuous functions  $f:\reals^n \to \reals^n$ whose inverse has a
continuous local selection; this follows from Brower's invariance of
domain theorem,  see the recent book \cite{KK} where
the interested reader can also find implicit function
theorems for nonsmooth mappings.

\bigskip 

\section{Aubin continuity and continuous local selections}

It is well known  that  $F$ is 
metrically regular at $\bx$ for $\by$
if and only if $F^{-1}$
has the so-called {\em Aubin property}
 at $\by$ for $\bx$:  
there exist $\kappa\in (0,\infty)$ 
together with  neighborhoods $U$ of $\bx$ and $V$ of $\by$ such that
$$
   F^{-1}(y')\cap U \,\subset\, F^{-1}(y) +\kappa\|y'-y\|\ball
        \text{   for all} y,y' \in V ;               \leqno(2.1)
$$
moreover the constants $\kappa$ in (1.3) and (2.1)
agree, that is, the  modulus $\,\reg F(\bx,\by)$ is also the infimum of
all $\kappa$ for which (2.1) holds.

The Aubin property (1.4) is a local property of a mapping
around a point in its graph, which is preserved after  a truncation
of the mapping  with a neighborhood of the reference point.
  However, such a truncation  may not be lower semicontinuous, in general.

{\bf Counterexample.}  Consider the mapping $A$ from $\reals$ to $\reals$
whose graph is the union of the graphs of the functions $x = y+ 1/k$
for $k = \pm 1,\pm 2, \cdots$ and the function $x=y$. 
This mapping is  Aubin continuous at zero for zero
 however, for any $\eps > 0$ and $\delta > 0$,
if we consider the restriction of the graph of $A$ in 
the box $[-\delta, \delta] \times [-\eps, \eps]$,
there will
be points in  this restriction
 with coordinates $(y,\eps)$ with $y < \delta$
that cannot be approached by a sequence $x_n$ so that $(y_n, x_n) \in \gph A$,
$y_n \to y$  and $\delta > y_n > y$.

In the following lemma we show 
that  if $A$ is  convex and closed valued locally around the reference point, then
 the mapping obtained by
 truncation of $A$ with a ball centered at $\bx$ with radius proportional
to the distance to $\by$
is lower semicontinuous in a neighborhood of $\by$. 

\proclaim Lemma 2.1. Consider a mapping $A:Y\tto X$ and any $(\by,\bx)\in \gph A$
and suppose that $A$ is Aubin continuous  at $\by$ for $\bx$
with a constant $\kappa$.  
 Let, for some $c > 0$,
the sets  $A(y)\cap \ball_c(\bx)$ be convex and closed for all $y \in \ball_c(\by)$.
The for any $\alpha > \kappa$
there exists $\beta > 0$ such that the mapping
$$ \ball_{\beta}(\by) \ni y \mapsto M_0(y):= \{
x \in A(y) \mid \|x - \bx\| \leq \alpha\|y - \by\| \} $$ 
 is nonempty, closed and convex valued, and 
 lower semicontinuous.
 
{\bf Proof.} Let $\kappa < \alpha $
and let $\ball_a(\bx)$ and $\ball_b(\by)$ 
be the neighborhoods of $\bx$ and $\by$, respectively,
that  are associated with the 
Aubin continuity of $A$ (metric regularity
of $A^{-1}$) with constant  $\kappa$. 
Choose $\beta > 0$
such that $\beta \leq a/\kappa$ that $\max \{a,\beta\} \leq  c$.
For such a $\beta$  the mapping
$M_0$ has  nonempty closed convex values. It remains to show that
$M_0$ is lower semicontinuous on $\ball_\beta(\by)$.

Let $(x,y) \in \gph M_0$ and $y_k \to y$, $y_k \in \ball_\beta(\by)$. 
First, let  $y = \by$.
Then $M_0(y) = \bx$ and from the Aubin continuity of $A$
there exists a sequence $x_k \in A(y_k)$ such that
$\|x_k - \bx\| \leq \kappa \|y_k - \by\|$. Thus $x_k \in M_0(y_k)$,
$x_k \to x$ as $k \to \infty$ and we are done in this case.

Now let $y \neq \by$. From the Aubin property
of $A$ there exists $\check x_k \in A(y_k)\cap \ball_c(\bx)$ such that
$$\|{\check x}_k - \bx \| \leq \kappa \|y_k - \by\|$$
and also there exists 
$\tilde x_k \in A(y_k)\cap \ball_c(\bx)$ such that
$$\|{\tilde  x}_k - x \| \leq \kappa \|y_k - y\|.$$
Let $$\mu_k = \frac{\|y_k - y\|}{(\alpha - \kappa)\|y_k - \by\| + \|y_k - y\|}.
$$
Then  $\mu_k \to 0$ as $k \to \infty$ and hence for large $k$ we have
$0 \leq \mu_k < 1$.
Let $x_k = \mu_k {\check x}_k + (1-\mu_k){\tilde x}_k$.
Then $x_k \in A(y_k)$. We have
\bas  \|x_k - \bx \| & \leq & \mu_k\|{\check x}_k - \bx\| + (1-\mu_k)\|{\tilde x}_k - \bx\| \\
 &   \leq & \mu_k \kappa \|y_k - \by\|+
(1-\mu_k)(\|{\tilde x}_k - x\| + \|x - \bx\|) \\
& \leq &  \mu_k \kappa \|y_k - \by\|+
(1-\mu_k)\kappa \|y_k - y\| + (1-\mu_k)\alpha \|y - \by\| \\
&  \leq &   \alpha \|y - \by\| - \mu_k (\alpha - \kappa) \|y_k - \by\|
+ (1-\mu_k)\kappa \|y_k - y\| 
 \leq \alpha \|y - \by\|,
\eas
because of the choice of $\mu_k$. Thus $x_k \in M_0(y_k)$ and 
since $x_k \to x$, 
the proof is complete. \qed

Lemma 2.1 allows us to apply the Michael selection theorem
to the mapping $M_0$ obtaining, in terms of a metrically regular mapping $F$,
the following theorem:

\proclaim Theorem 2.2. Consider a mapping $F:X\tto Y$ and
 $(\bx,\by)\in \gph F$ at which $\reg F(\bx\for\by) < \infty$. 
 Let, for some $c > 0$,
the sets  $F^{-1}(y)\cap \ball_c(\bx)$ be convex and closed for all $y \in \ball_c(\by)$.
Then the mapping $F^{-1}$ has a continuous local selection $x(\cdot)$
around $(\bx,\by)$ which is calm at $\by$
with $\calm x(\by) \leq  \reg F(\bx\for\by)$.

In the following section we will show that on the same assumptions, the conclusion of this theorem  holds when $F$ is perturbed by a function $G$ with a sufficiently small
Lipschitz constant.

\bigskip

\section{The inverse mapping theorem}

The metric regularity of  a  mapping $F$ is preserved when $F$ is 
perturbed by a function with a small Lipschitz constant.
This property of the metric regularity allows 
one to pass from linear to nonlinear and back and is 
usually identified with the Lyusternik-Graves theorem.
Here we refer to
 the following  version of this result proved in \cite{rad}.

\proclaim Theorem 3.1. 
Consider a mapping $F:X\tto Y$ with  $(\bx,\by) \in \gph F$ and let
$\,\gph F$ has a closed intersection with
a neighborhood of  $(\bx,\by)$.  Consider also a mapping $G:X\to Y$.  
If $\,\reg F(\bx\for\by) <\kappa <\infty$ and $\,\lip G(\bx)<\lambda < 
\kappa^{-1}$, then
$$ 
 \reg(F+G)(\bx\for\by+G(\bx)) <  (\kappa^{-1}-\lambda)^{-1}. 
$$

\noindent
 Recall that the  Lipschitz modulus $\lip G(\bx)$ of a 
single-valued mapping $G$ at a point $\bx$ is defined as 
$$
   \lip G(\bx) := \limsup_{x,x'\to \bx \atop x,x'\neq\bx}
                     {\|G(x')-G(x)\| \over \|x'-x\|}.    
$$
For a function $F:X\to Y$ which is strictly differentiable 
at $\bx$, Theorem 3.1 implies
$$
             \reg F(\bx\for F(\bx)) = \reg \nabla F(\bx).       
$$
 Since $\nabla F(\bx)$ is linear and bounded, from
the Banach open mapping theorem, $F$ is metrically regular at $\bx$ for 
$F(\bx)$ if and only if  $\nabla F(\bx)$ is surjective.  

The following theorem is the main result of
this paper. We show that
 if a mapping $F$ satisfies the condition of Theorem 2.2 in the previous section,
and hence has a continuous local selection around the reference point, then
for any function $G:X\to Y$ with $\,\lip G(\bx)<\reg F(\bx\for\by)$, 
the mapping $(F+G)^{-1}$  has a continuous and calm local selection around
the reference point which is also calm. Note that we claim the existence of 
a continuous selection of a possibly nonconvex valued mapping.
Also note that this result is an inverse function theorem for
a mapping $Q$ which can be represented as the sum $Q=F+G$ where
$F$ and $G$ have corresponding properties.
We prove this theorem using repeatedly the argument
in the proof of Lemma 2.1  in a way which resembles the
classical proofs of Lyusternik and Graves.

\proclaim Theorem 3.2. Consider a mapping $F:X\tto Y$ which
is metrically regular at $\bx$ for $\by$. Let
$\,\gph F$ have a closed intersection with
a neighborhood of  $(\bx,\by)$ and let for some $c > 0$
the mapping   
$\ball_c(\by) \ni y \mapsto F^{-1}\cap \ball_c(\bx)$ be convex valued.
Let $G:X\to Y$ satisfy  $\,\lip G(\bx)\reg F(\bx\for\by)< 1$. 
Then the mapping $(G+F)^{-1}$ has a continuous  local selection 
$x(\cdot)$ around $(\by+G(\bx), \bx)$ which is calm at $\by$ with 
$$  \calm x(\by) \leq    \frac{2\reg F(\bx\for\by)}
{1-\lip G(\bx) \reg F(\bx\for\by)} \, . \leqno{(3.1)}$$

\state Proof. Choose $\gamma$ that is greater than the right hand side of (3.1)
and let $\kappa,$ $ \alpha $ and
$\lambda$ be such that $\reg F(\bx\for \by)< \kappa < \alpha < 1/\lambda $, 
$\lambda > \lip G(\bx)$ and 
$2\kappa/(1-\alpha\lambda) \leq \gamma .$
For simplicity, we assume that
$G(\bx) = 0$. Let $\ball_a(\bx)$ and $\ball_b(\by)$ 
be the neighborhoods of $\bx$ and $\by$, respectively,
that  are associated with the 
metric regularity of $F$ at $\bx$ for $\by$ with  constant $\kappa$
and $G$ is Lipschitz continuous in $\ball_a(\bx)$
with constant $\lambda$.
Without loss of generality, let $\max\{a, b\} \leq c$.
Note that $F^{-1}(y) \neq \emptyset$ for any $y \in \ball_b(\by)$.
From Theorem~2.2, if we choose $b$ small enough, there exists a continuous
function $z_0:\ball_b(\bx) \to X$ such that
$$F(z_0(y))\ni y \Text{and} \|z_0(y) - \bx\| \leq \kappa \|y - \by\| $$
for all $y \in \ball_b(\by)$.  Choose a positive $\tau$ such that
$$ \tau \leq (1-\alpha\lambda) \min\bigg\{\frac{a}{2\kappa},\frac{b}{1+\kappa\lambda}\bigg\}.  \leqno{(3.2)}
$$
 Consider the mapping
$$ \ball_{\tau}(\by) \ni y \mapsto M_1(y):= \lset
x \in F^{-1}(y - G(z_0(y))) \mid \|x - z_0(y)\| \leq 
\kappa(1 + \kappa\lambda)  \|y - \by\| \rset.$$ 
The mapping $M_1$ is closed and convex valued and $(\by,\bx) \in \gph M_1$.
For $y \in \ball_\tau(\by)$ we have
$$\|y-G(z_0(y)) - \by\| \leq \tau + \lambda\|z_0(y) - \bx\| \leq
\tau + \lambda\kappa \tau \leq b ,$$
therefore the mapping $M_1$ is also nonempty valued. 
We will show that this mapping is lower semicontinuous in $\ball_\tau(\by)$.

Let $y \in \ball_{\tau}(\by)$ and $x \in M_1(y)$,
and let  $y_k \in \ball_\tau(\by)$, $y_k \to y$ as $k \to \infty$. 
First, assume that $y = \by$. Then $M_1(y) = \bx=z_0(\by)=z_0(y) \in F^{-1}(\by-G(\bx)) $ 
and from the Aubin
property of $F^{-1}$ there exists $x_k \in F^{-1}(y_k - G(z_0(y_k)))$
such that
$$ \|x_k - z_0(\by)\|  \leq  \kappa (\|y_k - \by\| + \|G(z_0(y_k)) - G(\bx)\|).
$$
Then, using the calmness of $z_0$ we obtain
 $$ \|x_k - z_0(y)\|  \leq  \kappa (\|y_k - \by\| + \lambda\|z_0(y_k) - \bx\|)
\leq  \kappa(1+\kappa\lambda)\|y_k - \by\|.$$
Thus $x_k \in M_1(y_k)$ and $x_k \to \bx$.

Now suppose that $y \neq \by$.
Using the calmness of $z_0$ and (3.2) we have 
$$\|y_k - G(z_0(y_k)) - \by \| \leq \|y_k - \by\|  + \|G(z_0(y_k)) - G(\bx)\|
\leq (1 +   \kappa\lambda)\tau < b \leqno{(3.3)}$$
 and
$$\|z_0(y_k) - \bx\| \leq \kappa \tau \leq a . \leqno{(3.4)}$$ 
Since $z_0(y_k) \in F^{-1}(y_k-G(\bx)))\cap \ball_a(\bx) $,
from the Aubin continuity of $F^{-1}$ there exists 
${\check x}_k \in F^{-1}(y_k-G(z_0(y_k)))$ such that
$$
\|\check{x}_k - z_0(y_k)\| \leq \kappa \| G(z_0(y_k))-G(\bx)\| \leq \kappa \lambda \|z_0(y_k)
- \bx\| 
\leq  \kappa^2 \lambda   \|y_k - \by\|. \leqno{(3.5)}
$$
Similarly, since 
$x \in F^{-1}(y -  G(z_0(y)))\cap \ball_a(\bx) $,
from the estimations (3.3) and (3.4) with $y_k$ replaced by $y$
and  from the Aubin continuity of $F^{-1}$ 
 there exists
$\tilde{x}_k \in F^{-1}(y_k - G(z_0(y_k))) $ such that
$$\|\tilde{x}_k - x\| \leq \kappa\|y_k - y\| + \kappa\lambda \|z_0(y_k) - z_0(y)\|, 
\leqno{(3.6)}$$
thus  $\tilde{x}_k\to x$ as $ k \to \infty$. Remembering that
 $y \neq \by$, for $k$ sufficiently large
we have $\|y_k-\by\| > \|y-\by\|/2 > 0$. 
For these large $k$, denote
$$\epsilon_k:=\frac{\|\tilde{x}_k - x\|+\|z_0(y_k) - z_0(y)\|
 }{\kappa\|y_k - \by\|-\kappa^2\lambda\|y_k - y\|}. $$
From (3.6)  and the continuity of $z_0$
we obtain that $\epsilon_k \to 0$. 
Let $\beta_k$   be a convergent to zero sequence such that
$\beta_k \geq \epsilon_k.$
From the local convexity of the values of $F^{-1}$, the point
$$x_k = \beta_k{\check x}_k + (1-\beta_k)\tilde{x}_k $$
is an element of $ F^{-1}(y_k-G(z_0(y_k)))$ and also
 $x_k \to x$ as $k \to \infty$.  Then (3.5), (3.6) and the choice of $\epsilon_k$
and $\beta_k$  yield
\begin{eqnarray*}
\|x_k - z_0(y_k)\| & \leq & \beta_k\|{\check x}_k - z_0(y_k)\|
+ (1-\beta_k)\|\tilde{x}_k - z_0(y_k)\| \\
& \leq & \beta_k \kappa^2\lambda\|y_k - y\|+ \beta_k \kappa^2\lambda\|y - \by\|\\
& &  + (1-\beta_k)(\|\tilde{x}_k - x\| +
\|x - z_0(y)\| + \|z_0(y) - z_0(y_k)\|)\\
& \leq & \beta_k \kappa^2\lambda\|y_k - y\|+ \beta_k \kappa^2\lambda\|y - \by\|\\
 & &+ (1-\beta_k)(\|\tilde{x}_k - x\|  
 + \kappa(1+\kappa\lambda)\|y - \by\| + \|z_0(y) - z_0(y_k)\|) \\
& \leq & \kappa(1+\kappa\lambda)\|y - \by\| - \beta_k\kappa\|y - \by\|
+ \beta_k \kappa^2\lambda\|y_k - y\| + \|\tilde{x}_k - x\| +  \|z_0(y) - z_0(y_k)\|\\
&=& \kappa(1+\kappa\lambda)\|y - \by\| -
(\beta_k-\epsilon_k)(\kappa\|y_k - \by\|-\kappa^2\lambda\|y_k - y\|) \\
&\leq &\kappa(1+\kappa\lambda)\|y - \by\| 
\end{eqnarray*}
for $k$ sufficiently large. We obtain that $x_k \in M_1(y_k)$ for such $k$. 
Hence  $M_1$ is lower semicontinuous in $\ball_\tau(\by)$.

The mapping $M_1$ is nonempty, closed and convex valued 
and lower semicontinuous  
in $\ball_\tau(\by)$, hence it
has a continuous selection  $z_1(\cdot):\ball_\tau(\by)
 \to X$ which  by  definition 
satisfies
$$ z_1(y) \in F^{-1}(y - G(z_0(y))) \Text{and}
\|z_1(y) - z_0(y)\| \leq \kappa(1+\lambda\kappa)\|y - \by\|. \leqno{(3.7)}$$
Hence   for all $y \in \ball_\tau(\by)$, 
$$ \|z_1(y) - \bx\| \leq \|z_1(y) - z_0(y)\|
+ \|z_0(y)-\bx\|  \leq \kappa(2+\kappa\lambda)\|y - \by\|. \leqno{(3.8)}$$
In particular, from (3.2)
$$\|z_1(y) - \bx\|  \leq a. $$
and also
$$\|y - G(z_1(y)) - \by\| \leq \|y - \by\| + \|G(z_1(y)) - G(\bx)\|
 \leq b. \leqno{(3.9)}$$

 Now consider the mapping
$$ \ball_{\tau}(\by) \ni y \mapsto M_2(y):= \lset
x \in F^{-1}(y - G(z_1(y))) \mid \|x - z_1(y)\| \leq 
\alpha\lambda\|z_1(y) - z_0(y)\| \rset.$$ 
Again, we will apply the Michael selection theorem to $M_2$ after proving that 
it is lower semicontinuous in $\ball_\beta(\by)$. 
Of course,  $M_2$ is nonempty, closed and convex valued and $(\by,\bx) \in \gph M_2$.

Let $y \in \ball_{\tau}(\by)$ and $x \in M_2(y)$,
and let  $y_k \in \ball_\tau(\by)$, $y_k \to y$ as $k \to \infty$. 
If $z_0(y) = z_1(y)$ then $M_2(y) = \{z_1(y)\}$ and therefore $x = z_1(y)$,
and since 
$$z_1(y_k) \in F^{-1}(y_k - G(z_0(y_k)))\cap B_a(\bx) \Text{ and }
y_k - G(z_0(y_k)) \in \ball_b(\by)$$ 
 the Aubin  continuity
 of $F^{-1}$ implies that there exists $x_k \in F^{-1}(y_k - G(z_1(y_k)))$
such that
$$ \|x_k - z_1(y_k)\|  \leq  \kappa \|G(z_1(y_k)) - G(z_0(y_k))\| 
\leq \alpha\lambda \|z_1(y_k) - z_0(y_k)\|. $$
Hence $x_k \in M_2(y_k)$ and from the continuity
of the functions $z_0$ and $z_1$ we obtain that
$x_k \to z_1(y) =x$, thus $M_2$ is lower semicontinuous.

Now let $z_0(y) \neq z_1(y)$.
The estimations (3.8) and (3.9) clearly hold for $y$ replaced by $y_k$, and
 since $z_1(y_k) \in F^{-1}(y_k-G(z_0(y_k)))\cap \ball_a(\bx) $,
the Aubin continuity of $F^{-1}$ implies the existence of
${\check x}_k \in F^{-1}(y_k-G(z_1(y_k)))$ such that
$$
\|\check{x}_k - z_1(y_k)\| \leq \kappa \|G(z_1(y_k)) - G(z_0(y_k))\| \leq \kappa \lambda 
\|z_1(y_k) - z_0(y_k)\|. \leqno{(3.10)}$$
Also,  taking into account (3.8)  and (3.9) and 
the inclusion $x \in F^{-1}(y -  G(z_1(y)))\cap \ball_a(\bx) $,
the Aubin continuity of $F^{-1}$ yields that
there exists
$\tilde{x}_k \in F^{-1}(y_k - G(z_1(y_k))) $ such that
$$\|\tilde{x}_k - x\| \leq \kappa(\|y_k - y\| + \lambda \|z_1(y_k) - z_1(y)\|)  \to 0
\text{  as  } k \to \infty . \leqno{(3.11)} $$

Let $\beta_k $ be an arbitrary sequence of positive numbers
that is convergent to zero as $k \to \infty$ and let
$$\epsilon_k:=\frac{\beta_k\kappa\lambda\|z_1(y_k) - z_0(y_k)\| + (1-\beta_k)
(\|\tilde{x}_k - x\| + \|z_1(y_k)- z_1(y)\|) }{ \lambda \|z_1(y_k) - z_0(y_k)\|} .$$
Note that $\|z_1(y_k) - z_0(y_k)\| \geq \|z_1(y) - z_0(y)\|/2>0$ for all
large $k$ and therefore $\epsilon_k \to 0$ as $k \to \infty$. Let
$$x_k = \beta_k{\check x}_k + (1-\beta_k)\tilde{x}_k .$$
Since $F^{-1}$ is locally convex valued, we have $x_k \in F^{-1}(y_k-G(z_1(y_k)))$ and,
since $\tilde{x}_k \to x$ and $\beta_k \to 0$,
we obtain  $x_k \to x$ as $k \to \infty$.  From (3.10), (3.11)
and the choice of $\epsilon_k$
and $\beta_k$ we have
\begin{eqnarray*}
\|x_k - z_1(y_k)\| & \leq & \beta_k\|{\check x}_k - z_1(y_k)\|
+ (1-\beta_k)\|\tilde{x}_k - z_1(y_k)\| \\
& \leq & \beta_k \kappa\lambda \|z_1(y_k) - z_0(y_k)\| +
(1-\beta_k)(\|\tilde{x}_k - x\| + \|x - z_1(y)\| + \|z_1(y) - z_1(y_k)\|) \\
& \leq &  \beta_k \kappa\lambda \|z_1(y_k) - z_0(y_k)\|  \\ & & +
(1-\beta_k)(\|\tilde{x}_k - x\| + \kappa\lambda\|z_0(y) - z_1(y)\| + \|z_1(y) - z_1(y_k)\|) \\
& \leq & \alpha\lambda\|z_0(y) - z_1(y)\| - (\alpha - \kappa)\lambda\|z_0(y) - z_1(y)\| \\
& & \,\, + \beta_k\kappa\lambda\|z_1(y_k) - z_1(y)\| + 
(1 -\beta_k)(\|z_1(y_k) - z_1(y)\| + \|\tilde{x}_k - x\|)\\
& \leq & \alpha\lambda\|z_0(y) - z_1(y)\| - (\alpha - 
\kappa - \epsilon_k)\lambda\|z_0(y) - z_1(y)\| \\
& \leq &  \alpha\lambda\|z_0(y) - z_1(y)\| 
\end{eqnarray*}
for sufficiently large $k$. 
We obtain that $x_k \in M_2(y_k)$ for all large $k$
and since $x_k \to x$, the mapping  $M_2$ is lower semicontinuous in $\ball_\tau(\by)$.
Hence it
has a continuous selection  $z_2(\cdot):\ball_\tau(\by) \to X$; by definition it 
satisfies
$$ z_2(y) \in F^{-1}(y - G(z_1(y))) \Text{and}
\|z_2(y) - z_1(y)\| \leq \alpha\lambda\|z_1(y) - z_0(y)\| 
\text{for all} y \in \ball_\tau(\by).  $$

The induction step is  analogous. 
Let $z_0, z_1$ and $z_2$ be as above and
suppose we have also found  functions $  z_3, z_4,\cdots, z_n$, 
 such that any $z_j$, 
$j =3,4, \cdots, n$, 
is a continuous local selection of the  mapping
$$\ball_{\tau}(\by) \ni y \mapsto M_{j}(y):= \lset
x \in F^{-1}(y - G(z_{{j-1}}(y)))
 \mid \|x - z_{j-1}(y)\| \leq \alpha\lambda \|z_{j-1}(y)- z_{j-2}(y)\|\rset $$
Then for $y \in \ball_\tau(\by)$ we obtain
$$\|z_j(y) - z_{j-1}(y)\| \leq (\alpha\lambda)^{j-1}\|z_1(y) - z_0(y)\| $$
and therefore, using Theorem 2.2 and (3.7),
$$
\|z_j(y) - \bx\| \leq 
\sum_{i=1}^{j}(\alpha \lambda)^{i-1} 
\|z_1(y) - z_0(y)\| + \|z_0(y) - \bx\|
 \leq \frac{2\kappa}{1-\alpha \lambda}\|y - \by\|. $$
Hence, from (3.2),  for
$j =3,4 \cdots, n$, 
$$
\|z_j(y) - \bx\| \leq \frac{2\kappa\tau}{1-\alpha \lambda} \leq a \leqno{(3.12)}
$$
and
$$\|y - G(z_j(y)) - \by\| \leq \tau + \|z_j(y) - \bx\|
\leq \tau +  \frac{2\kappa\lambda\tau}{1-\alpha \lambda} \leq b .  \leqno{(3.13)}$$
Consider the mapping
$$ \ball_{\tau}(\by) \ni y \mapsto M_{n+1}(y):= \lset
x \in F^{-1}(y - G(z_n(y))) \mid \|x - z_n(y)\| \leq 
\alpha\lambda\|z_n(y) - z_{n-1}(y)\| \rset.$$ 
which is nonempty,  closed and convex valued and $(\by,\bx)$ is in its graph.
Let $y \in \ball_{\tau}(\by)$ and $x \in M_{n+1}(y)$,
and let  $y_k \in \ball_\tau(\by)$, $y_k \to y$ as $k \to \infty$. 
If $z_{n-1}(y) = z_n(y)$ then $M_{n+1}(y) = \{z_n(y)\}$ and hence $x = z_n(y)$,
and from $z_n(y_k) \in F^{-1}(y_k - G(z_{n-1}(y_k)))\cap B_a(\bx)$
and $y_k - G(z_{n-1}(y_k)) \in \ball_b(\by)$,
and  using the Aubin
property of $F^{-1}$, we obtain 
 that there exists $x_k \in F^{-1}(y_k - G(z_n(y_k)))$
such that
$$ \|x_k - z_n(y_k)\|  \leq  \kappa \|G(z_n(y_k)) - G(z_{n-1}(y_k))\| 
\leq \alpha\lambda \|z_n(y_k) - z_{n-1}(y_k)\|. $$
Therefore $x_k \in M_{n+1}(y_k)$,
$x_k \to z_1(y) =x$ as $k \to \infty$, and hence $M_2$ is lower semicontinuous.

Let $z_n(y) \neq z_{n-1}(y)$.
From (3.12) and (3.13) for $y = y_k$,
 since $z_n(y_k) \in F^{-1}(y_k-G(z_{n-1}(y_k)))\cap \ball_a(\bx) $,
the Aubin continuity of $F^{-1}$ implies the existence of
${\check x}_k \in F^{-1}(y_k-G(z_n(y_k)))$ such that
$$
\|\check{x}_k - z_1(y_k)\| \leq \kappa \|G(z_n(y_k)) - G(z_{n-1}(y_k))\| \leq \kappa \lambda 
\|z_n(y_k) - z_{n-1}(y_k)\|. $$
Similarly, the estimations (3.12)  and (3.13) 
and   the Aubin continuity of $F^{-1}$ yield that,
since 
$x \in F^{-1}(y -  G(z_n(y)))\cap \ball_a(\bx) $, there exists
$\tilde{x}_k \in F^{-1}(y_k - G(z_n(y_k))) $ such that
$$\|\tilde{x}_k - x\| \leq \kappa(\|y_k - y\| + \lambda \|z_n(y_k) - z_n(y)\|)
\leq \kappa(\|y_k - y\| + \|G(z_n(y_k)) - G(z_n(y))\|)  \to 0
\text{  as  } k \to \infty . $$
Choose an arbitrary sequence  $\beta_k \to 0$ 
which is convergent to zero sequence and let
$$\epsilon_k:=\frac{\beta_k\kappa\lambda\|z_n(y_k) - z_{n-1}(y_k)\| + (1-\beta_k)
(\|\tilde{x}_k - x\| + \|z_n(y_k)- z_n(y)\|) }{ \lambda \|z_n(y_k) - z_{n-1}(y_k)\|} .$$
Then $\epsilon_k \to 0$ as $k \to \infty$. Taking
$$x_k = \beta_k{\check x}_k + (1-\beta_k)\tilde{x}_k, $$
we estimate the distance 
$\|x_k - z_n(y_k)\|$ in the same way as in the first step,
by just replacing  $z_1$ by $z_n$ and
$z_0$ by $z_{n-1}$; then  we conclude that $x_k \in M_{n+1}(y_k)$ for all large $k$.
  Since   $x_k \to x$ as $k \to \infty$,
the mapping $M_{n+1}$ is lower semicontinuous in $\ball_\tau(\by)$.
Hence $M_{n+1}$
has a continuous selection  $z_{n+1}(\cdot):\ball_\tau(\by) \to X$ which
then satisfies
$$ z_{n+1}(y) \in F^{-1}(y - G(z_n(y))) \Text{and}
\|z_{n+1}(y) - z_n(y)\| \leq \alpha\lambda\|z_n(y) - z_{n-1}(y)\| . $$
Putting all  together we have
$$\|z_{n+1}(y) - z_n(y)\| \leq (\alpha\lambda)^n\|z_1(y) - z_{0}(y)\|  \leqno{(3.14)}$$
and
$$\|z_{n+1}(y) - \bx \| \leq 
\sum_{i=1}^{n+1}(\alpha \lambda)^{i-1} 
\|z_1(y) - z_0(y)\| + \|z_0(y) - \bx\|
 \leq \frac{2\kappa}{1-\alpha \lambda}\|y - \by\|, $$
and the induction step is complete. 

Thus, we obtain an infinite sequence of functions $z_0, z_1, \cdots, z_n, \cdots$ for which
(3.7) and (3.14) yield
$$\sup_{y \in \ball_\tau(\by)}\|z_{n+1}(y) - z_n(y)\| 
\leq (\alpha\lambda)^n \kappa(1+\lambda\kappa)\tau , $$
therefore, since $\alpha\lambda < 1$,
$\{z_n\}$ is a Cauchy sequence in the space of functions that are
continuous on $\ball_{\tau}(\by)$ equipped with the supremum norm. Then
this sequence  has  a limit $x(\cdot)$ which is a continuous
function in $\ball_{\tau}(\by)$ and  satisfies 
$$ x(y) \in F^{-1}(y - G(x(y)))  \text{and}  
\|x(y) - \bx\| \leq 
\frac{2\kappa}{1-\alpha\lambda}\|y - \by\|
\leq \gamma \|y - \by\|  
$$ 
for all  $y \in \ball_{\tau}(\by)$. This completes the proof. \qed

{\bf Proof of Theorem 1.2.} 
Apply Theorem 3.2 with $F(x)  = \nabla f(\bx)(x - \bx)$ and
$G(x) =
f(x) -  \nabla f(\bx)(x - \bx)$. 
Metric regularity of $F$ is equivalent to the surjectivity of
$\nabla f(\bx)$ and, since  $\nabla f(\bx)$ is linear and continuous,
the mapping  $F^{-1}$ is convex and closed valued. The mapping $G$
 has $\lip G(\bx) = 0$ and
finally $F+G =  f$. \qed

In this case a formula for the modulus of metric regularity
of the linear and bounded mapping $\nabla f(\bx)$ is available
from \cite{rad}, Example 1.1, and $\lip G(\bx) = 0$; then
the upper bound for the modulus of calmness has the form
$$\calm x(\by) \leq 2\sup \{ d(0, \nabla f(\bx)^{-1}(y)) \mid y \in \ball \}.$$

\bigskip

\section{Applications}

Theorems  3.2  can be also stated in a corresponding ``implicit function" form
as follows:

\proclaim Theorem 4.1. Let $X, Y$ and $Z$
be Banach spaces. Consider a mapping $F:X\tto Y$ and $(\bx,\by)\in \gph F$ 
which satisfies the conditions in Theorem 3.2. 
Consider also a mapping $G:X\times Z\to Y$ 
which is continuous in a neighborhood of $(\bx,\bp)$ and
with  $\,\lip_x G(\bx,\bp)<{\reg F(\bx\for\by)}^{-1}$
(here the Lipschitz modulus of $G(x,p)$ is with respect to $x$ 
where $\limsup$ is also with respect to
$p \to \bp$). Then there exists a neighborhood $U$ of $\bx$ and $ P$ of $\bp$, 
 a continuous function $x(\cdot): P \to U$, and a constant $\gamma$ such that
$$ 
  \by \in G(x(p), p) + F(x(p)) \, \Text{ and } 
  \|x(p) - \bx\| = \gamma \|G(\bx, p) -G(\bx,\bp)\| \text{  for every } p \in  P.
$$

{\bf Sketch of proof.} The proof is parallel to the proof of Theorem 3.2.
First we choose $\kappa,$ $ \alpha $ and
$\lambda$ such that $\reg F(\bx\for \by)< \kappa < \alpha < 1/\lambda $ and 
$\lambda > \lip_x G(\bx,\bp)$ and neighborhoods of $\bx$,  $\by$
and $\bp$ 
that  are associated with the 
metric regularity of $F$ at $\bx$ for $\by$ with  constant $\kappa$
and $G$ is Lipschitz continuous with respect to $x$ 
with constant $\lambda$ uniformly in $p$. For simplicity, let $G(\bx, \bp) = 0$
and $\by = 0$.
By appropriately choosing a sufficiently small radius $\tau$
of a ball around $\bp$,
 we construct an infinite sequence of continuous functions
$z_j: \ball_\tau(\bp) \to X$, $j = 0,1, \cdots,$ that is uniformly in
$\ball_\tau(\bp) $ convergent to a function $x(\cdot)$ that satisfies
the conclusion of the theorem.
The first $z_0$ is obtained with the help of Theorem 2.2 and 
satisfies
$$z_0(p) \in F^{-1}(-G(\bx, p)) \Text{ and } \|z_0(p) - \bx\|
\leq \kappa \|G(\bx, p)\| .$$
The function $z_1$ is a continuous selection of the mapping
$$\ball_{\tau}(\bp) \ni y \mapsto M_{1}(p):= \lset
x \in F^{-1}(- G(z_0(p), p))
 \mid \|x - z_0(p)\| \leq \kappa(1+\lambda\kappa)\|G(\bx, p)\|\rset $$
while, analogously, $z_j$ is a continuous selection of 
$$\ball_{\tau}(\bp) \ni y \mapsto M_{j}(p):= \lset
x \in F^{-1}(- G(z_{{j-1}}(p),p))
 \mid \|x - z_{j-1}(p)\| \leq \alpha\lambda \|z_{j-1}(p)- z_{j-2}(p)\|\rset .$$
Then for $p \in \ball_\tau(\bp)$ we obtain
$$\|z_j(p) - z_{j-1}(p)\| \leq (\alpha\lambda)^{j-1}\|z_1(p) - z_0(p)\| $$
and also, for an appropriate $\gamma > 0$, 
$$
\|z_j(y) - \bx\| \leq 
\gamma\|G(\bx, p)\|. $$
Passing to the limit  completes the proof. \qed

If a mapping $F:X \tto Y$ has convex and closed graph, then, by the Robinson-Ursescu
theorem, the metric regularity of $F$ at $\bx$ for $\by$ is equivalent to
the condition $\by \in \Int {\rm Im} F$. For such mapping we obtain
the following corollary of Theorem 3.2:

\proclaim Corollary 4.2. Let $F:X \tto Y$ have convex and closed graph, let $(\by,\by)
\in \gph F$
and let $f:X\to Y$ be strictly differentiable
at $\bx$. Let the strict derivative $\nabla f(\bx)$ together with $F$
satisfy the condition 
$$\by \in \Int  {\rm Im}(f(\bx) + \nabla f(\bx)(\cdot - \bx)
+ F(\cdot)). \leqno{(4.1)}$$
Then there exist neighborhoods $U$ of $\bx$ and $V$ of $\by$, 
 a continuous function $x(\cdot):V \to U$, and a constant
$\gamma$ such that
$$(f + F)(x(y)) \ni y \Text{ and } \|x(y) - \bx\| \leq \gamma \|y - \by\| 
\text{  for every } y \in  V.$$ 

An implicit function version if the above corollary easily follows from 
Theorem 4.3.

As a more specific  application we consider the following 
controlled boundary value problem:
$$\dot{x}(t) = f(x(t), u(t)), \quad x(0) = 0, \,\, x(1) = b, \leqno{(4.2)}$$
where $f:\reals^n \times \reals^m \to \reals^n$ is a smooth function,
the control $u(t) \in {\cal U}$ where ${\cal U}$ is convex and  compact subset of $\reals^m$.
The pair $(x,u)$ is a {\em feasible } solution of (4.2) when it satisfies 
the differential equation and $u(t) \in {\cal U}$ for almost every $
t \in [0,1]$,  and also $x \in  
W^{1,\infty}_0([0,1], \reals^n)$, the space of all Lipschitz continuous functions $x$
with values in $\reals^n$ and 
with $x(0) = 0$, and $u \in L^\infty([0,1], \reals^m)$, the space of all essentially bounded and measurable functions with values in $\reals^m$. 
We equip $L^\infty$ with the esssup norm and
 and $W^{1,\infty}$ with the norm
$\|x\|_{1,\infty} = \|\dot x\|_\infty$. 
For simplicity, we assume that $f(0,0) = 0$ and $0 \in {\cal U}$ and take
$(0,0)$ as the reference solution. 

We apply Corollary 4.2 with the following specifications: 
$X = W^{1,\infty}_0([0,1], \reals^n)\times L^\infty([0,1], \reals^m)$
and  $Y = L^\infty([0,1], \reals^n)\times\reals^n,$
$F(x,u) = (Ax+Bu - \dot{x}, x(1))$ where $A = \nabla_xf(0,0), B = \nabla_uf(0,0)$,
$G(x,u) = (f(x,u) -Ax-Bu, 0).$ Then $(G+F)(x,u)= (f(x,u)- \dot{x}, x(1)) $.
Clearly, $F$ has convex and closed graph.
The condition (4.1) is equivalent to the
following: there exists an $\epsilon > 0$ such that for any $(y,b)$,
$y \in L^{\infty}([0,1], \reals^n) $ and $b \in \reals^n$ with
$\|y\|_\infty + \|b\| < \epsilon $, there exists a feasible solution
$(x,u)$ of the linearized boundary value problem
$$\dot{x}(t)= Ax(t) + Bu(t) - y(t), \quad x(0) = 0, x(1) = b. $$ 
The latter condition in turn is  equivalent to the existence
of a feasible solution of 
$$\dot{x}(t)= Ax(t) + Bu(t), \quad x(0) = 0,\, x(1) = b. $$
for all $b$ with sufficiently small norm. 
This property of the linear system is so-called 
{\em null-controllability} and can be equivalently written as
$$ 0 \in \Int\int_0^1e^{At}B{\cal U}dt,$$
where the integral is in the sense of Aumann. If $0 \in \Int {\cal U}$,
the null-controllability is equivalent to the
rank condition ${\rm rank}[B, AB, \cdots, A^{n-1}B] = n$.

Summarizing, of the linearization $\dot{x}(t)= Ax(t) + Bu(t)$ of (4.2)
is null-controllable with controls from ${\cal U}$, then 
there is a continuous function $b \mapsto (x,u)$ from a neighborhood $V$
of zero in $\reals^n$ to the product $W^{1,\infty}_0([0,1], \reals^n)
\times L^{\infty}([0,1], \reals^m)$ such that for each $b \in V$,
$(x, u)$ is a solution of the controlled boundary value problem (4.2),
moreover $(x,u)$ is calm at zero.
Note that neither Theorem 1.1 nor Theorem 1.2
may  be applied to this problem.

\enddocument